\documentclass[10pt]{article}

\usepackage{latexsym}
\usepackage{bezier}

\begin{document}

%

\newcommand{\EE}{\mathop{\rm I\! E}\nolimits}
\newcommand{\E}{\mathop{\rm E}\nolimits}
\newcommand{\I}{\mathop{\rm Im\, }\nolimits}
\newcommand{\Str}{\mathop{\rm Str\, }\nolimits}
\newcommand{\Sdet}{\mathop{\rm Sdet\, }\nolimits}
\newcommand{\STr}{\mathop{\rm STr\, }\nolimits}
\newcommand{\R}{\mathop{\rm Re\, }\nolimits}
\newcommand{\CC}{\mathop{\rm C\!\!\! I}\nolimits}
\newcommand{\FF}{\mathop{\rm I\! F}\nolimits}
\newcommand{\KK}{\mathop{\rm I\! K}\nolimits}
\newcommand{\LL}{\mathop{\rm I\! L}\nolimits}
\newcommand{\MM}{\mathop{\rm I\! M}\nolimits}
\newcommand{\NN}{\mathop{\rm I\! N}\nolimits}
\newcommand{\PP}{\mathop{\rm I\! P}\nolimits}
\newcommand{\QQ}{\mathop{\rm I\! Q}\nolimits}
\newcommand{\RR}{\mathop{\rm I\! R}\nolimits}
\newcommand{\ZZ}{\mathop{\sl Z\!\!Z}\nolimits}
\newcommand{\integer}{\mathop{\rm int}\nolimits}
\newcommand{\erf}{\mathop{\rm erf}\nolimits}
\newcommand{\diag}{\mathop{\rm diag}\nolimits}
\newcommand{\fl}{\mathop{\rm fl}\nolimits}
\newcommand{\eps}{\mathop{\rm eps}\nolimits}
\newcommand{\var}{\mathop{\rm var}\nolimits}


\newcommand{\pfeil}{\rightarrow}

\newcommand{\kat}{{\cal C}}
\newcommand{\rmat}{{\cal R}}
\newcommand{\oalg}{{\cal A}}
\newcommand{\falg}{{\cal F}}
\newcommand{\eich}{{\cal G}}
\newcommand{\hilb}{{\cal H}}
\newcommand{\calm}{{\cal M}}
\newcommand{\mod}{{\cal M}}
\newcommand{\kegel}{{\cal O}}
\newcommand{\kegels}{{\cal K}}
\newcommand{\bigrho}{\rho_\oplus}
\newcommand{\bigphi}{\phi_\oplus}

\newcommand{\tprod}{\otimes}

\newcommand{\horab}{\rule[-1mm]{0pt}{5mm}}
\newcommand{\iso}{\stackrel{\sim}{=}}
\newcommand{\quer}[1]{\overline{#1}}
\newcommand{\schlange}[1]{\widetilde{#1}}
\newcommand{\CVO}[3]{{{#3 \choose #1\hspace{3pt}#2}}}
\newcommand{\clebsch}[6]{{\left[{#1\atop#4}\hspace{3pt}{#2\atop#5}
                                     \hspace{3pt}{#3\atop#6}\right]}}
\newcommand{\sjsymbol}[6]{{\left\{{#1\atop#4}\hspace{3pt}{#2\atop#5}
                                     \hspace{3pt}{#3\atop#6}\right\}}}
\newcommand{\spann}{{{\rm span}}}
\newcommand{\Nat}{{{\rm Nat}}}
\newcommand{\Mor}{{{\rm Mor}}}
\newcommand{\End}{{{\rm End}}}
\newcommand{\ev}{{{\rm ev}}}
\newcommand{\coev}{{{\rm coev}}}
\newcommand{\id}{{{\rm id}}}
\newcommand{\Id}{{{\rm Id}}}
\newcommand{\Vec}{{{\rm Vec}}}
\newcommand{\Rep}{{{\rm Rep}}}
\newcommand{\ZB}{{{\rm ZB}}}
\newcommand{\TB}{{{\rm TB}}}
\newcommand{\HB}{{{\rm HB}}}
\newcommand{\BB}{{{\rm BB}}}
\newcommand{\ket}[1]{|#1\rangle}
\newcommand{\vak}{\ket{0}}
\newcommand{\bild}[3]{{        
  \unitlength1mm
  \begin{figure}[ht]
  \begin{picture}(120,#1)\end{picture}
  \caption{\label{#3}#2}
  \end{figure}
}}

\newcommand{\dottedline}[5]{
                   \multiput(#1,#2)(#3,#4){#5}{\circle*{0.5}}
                        }

\newenvironment{bew}{\textbf{Proof:}}{\hfill$\Box$}

\newtheorem{bem}{\textit{Remark}}
\newtheorem{bsp}{\textit{Example} }
\newtheorem{axiom}{\textit{Axiom}  }
\newtheorem{de}{\textit{Definition}	}
\newtheorem{satz}{\textit{Proposition} }
\newtheorem{lemma}[satz]{\textit{Lemma}	}
\newtheorem{kor}[satz]{\textit{Corollary}}
\newtheorem{theo}[satz]{\textit{Theorem}   }

\newcommand{\sbegin}[1]{\small\begin{#1}}
\newcommand{\send}[1]{\end{#1}\normalsize}

\sloppy

\title{The Potts Model with a Reflecting Boundary}
\author{Reinhard H\"aring-Oldenburg\\
Mathematisches Institut \\ Bunsenstr. 3-5\\
 37073 G\"ottingen, Germany\\
email: haering@cfgauss.uni-math.gwdg.de}
\date{February 6, 1996}
\maketitle

\begin{abstract}
A Potts model with a reflecting boundary is introduced and
it is shown that its  partition function can be
expressed as a Markov trace on the Temperley-Lieb Algebra of
Coxeter type B.
\end{abstract}

\setlength{\unitlength}{1mm}

\section{Introduction}

Statistical and quantum field theoretic models with boundary
conditions have recently attracted increasing interest. We believe that
the theory of knots and braids of Coxeter type B plays an
equally important role in these situations as the ordinary knot
theory does in models without reflecting boundaries.

\begin{figure}[ht]
\begin{picture}(150,40)
\linethickness{0.5mm}

\put(20,10){\circle*{2}}
\put(20,20){\circle*{2}}
\put(20,30){\circle*{2}}
\put(30,10){\circle*{2}}
\put(30,20){\circle*{2}}
\put(30,30){\circle*{2}}
\put(40,10){\circle*{2}}
\put(40,20){\circle*{2}}
\put(40,30){\circle*{2}}
\put(50,10){\circle*{2}}
\put(50,20){\circle*{2}}
\put(50,30){\circle*{2}}

\put(20,10){\line(0,1){20}}
\put(30,10){\line(0,1){20}}
\put(40,10){\line(0,1){20}}
\put(50,10){\line(0,1){20}}

\put(20,10){\line(1,0){30}}
\dottedline{10}{10}{1}{0}{10}
\put(20,20){\line(1,0){30}}
\dottedline{10}{20}{1}{0}{10}
\put(20,30){\line(1,0){30}}
\dottedline{10}{30}{1}{0}{10}

\multiput(5,4)(0,3){10}{\line(1,1){5}}

\linethickness{1mm}
\put(10,5){\line(0,1){30}}
\put(32,2){\mbox{(a)}}

\put(102,2){\mbox{(b)}}
\linethickness{0.5mm}
\multiput(75,4)(0,3){10}{\line(1,1){5}}
\linethickness{1mm}
\put(80,5){\line(0,1){30}}
\linethickness{0.5mm}
\put(80,20){\circle*{2}}
\put(80,30){\circle*{2}}
\put(90,10){\circle*{2}}
\put(100,10){\circle*{2}}
\put(100,20){\circle*{2}}
\dottedline{80}{10}{1}{0}{10}
\put(80,25){\oval(10,10)[r]}
\put(80,20){\line(1,-1){10}}
\put(90,10){\line(1,1){10}}
\put(100,10){\line(0,1){10}}
\put(90,10){\line(1,0){10}}

\end{picture}
\caption{\label{lattice} (a) A lattice with boundary (b) A boundary graph}
\end{figure}
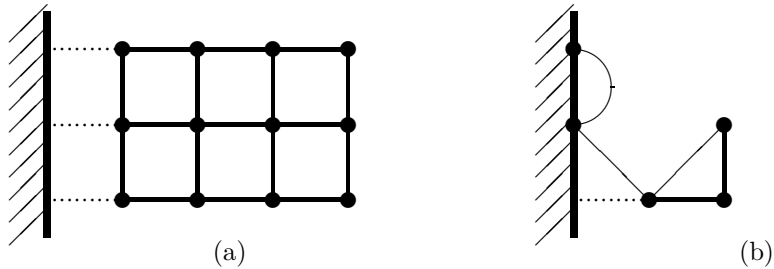

In this paper we investigate  a generalisation of the ordinary
Potts  model \cite{kauff} by including a  reflecting boundary.
Hence we have besides the usual lattice of sites a wall
interacting with it.
We visualise this as in Figure \ref{lattice}a. The dotted lines
indicate interaction bonds with the wall while the solid lines are
usual bonds between sites. Each site supports a 'spin' which
may occupy one of $f$ states.

\section{The Potts model on an arbitrary graph}\label{pottsmod}

It is useful to consider the lattice as a special kind of graph and
consider Potts models defined on arbitrary graphs. We further need to 
allow some vertices to be located on the wall. Such configurations
occur as intermediate states in the recursive
procedure for calculating the partition function we are going to present.

A Potts model with boundary
lives on a graph $G=(V,B)$ with a set of vertices $V$ and a
set of edges (bonds) $B$. Some of the vertices are defined to lie
on the boundary. They are bound with special bonds to the
reflecting wall. Formally we write the set of bonds
$B=B_1\cup B_0$ as a
disjoint union of inner bonds $B_1\subset V\times V$ 
and boundary bonds $B_0\subset V$.
We also single out a subset $V_0\subset V$ of vertices that lie on
the wall. These vertices shall have no degree of freedom: 
The associated spin is always in the ground state.

The whole structure shall be called a {\rm\bf graph with boundary}.

In figure \ref{lattice} dotted lines (resp solid lines)  indicate
 boundary bonds (resp. internal bonds). Note that bonds ending at a vertex
on the wall are internal bonds. 

We allow each vertex in $V\backslash V_0$ to aquire one of $f$ states. 
Vertices in $V_0$ are fixed to occupy state $0$. Hence
the set of all states is 
${\cal S}=\{S:V\rightarrow\{0,1,\ldots,f-1\}\mid S(i)=0, i\in V_0\}$.
The partition function (with $k$ being the Boltzmann constant and $T$
the temperature)   is
\begin{equation}\label{zdef}
Z_G=\sum_{S\in{\cal S}}{\rm exp}\left(\frac{-E(S)}{kT}\right)
\end{equation}
with the Hamiltonian  ($\delta(x,y)$ is the Kronecker symbol
with values $0,1$)
\begin{equation}
E(S)=\sum_{(i,j)\in B_1}\delta(S_i,S_j)+
\kappa\sum_{i\in B_0}(1-\delta(0,S_i))
\end{equation}

The first term in $E(S)$ is the usual Hamiltonian of the Potts model.
The second term introduces the boundary condition.

\begin{eqnarray}
Z_G&=&\sum_{S\in{\cal S}}{\rm exp}\left(-\frac{1}{kT}
  \sum_{(i,j)\in B_1}\delta(S_i,S_j)-\frac{\kappa}{kT}
     \sum_{i\in B_0}1-\delta(0,S_i)\right)\nonumber\\
&=&\sum_{S\in{\cal S}}\prod_{(i,j)\in B_1}{\rm exp}\left(-\frac{1}{kT}
  \delta(S_i,S_j)\right)\prod_{i\in B_0}{\rm exp}\left(\frac{\kappa}{kT}
     (\delta(0,S_i)-1)\right)\nonumber\\
&=&\sum_{S\in{\cal S}}\prod_{(i,j)\in B_1}\left(
 1+\delta(S_i,S_j)(e^{-\frac{1}{kT}}-1)\right)
  \prod_{i\in B_0}\left(e^{-\frac{\kappa}{kT}}+
     \delta(0,S_i)(1-e^{-\frac{\kappa}{kT}})\right)\label{z1}
\end{eqnarray}

We introduce the following short cuts 
\begin{equation}
A:=1\qquad B:=e^{-\frac{1}{kT}}-1\qquad C:=e^{-\frac{\kappa}{kT}}
\qquad D:=1-C
\end{equation}
and take a factor $C$ out of the sum
\begin{eqnarray}
Z_G&=&\sum_{S\in{\cal S}}\prod_{(i,j)\in B_1}\left(
 1+\delta(S_i,S_j)B\right)
  \prod_{i\in B_0}\left(C+
     \delta(0,S_i)D\right)\label{zzz}\\
&=&C^{\#B_0}
\sum_{S\in{\cal S}}\prod_{(i,j)\in B_1}\left(
 1+\delta(S_i,S_j)B\right)
  \prod_{i\in B_0}\left(1+
     \delta(0,S_i)D/C\right)
\end{eqnarray}

From this expression it can easily be read off how $Z_G$  changes 
when $G$ is changed.

Lets fix a single inner bond and let $G_d$ be the
graph obtained from $G$ by deleting this bond. Further let
$G_c$ be the graph $G$ where the bond is contracted to a point so that
two vertices merge into one and the bond disappears.
Then we have from (\ref{zzz}):
\begin{equation} Z_G=Z_{G_d}+BZ_{G_c}\label{g1}
\end{equation}
It is irrelevant whether any of the two vertices of the bond 
lie on the wall. 

A vertex $i\in V\backslash V_0$ without any bonds may freely occupy all states
without contributing to the energy. So if $G_1$ is obtained
from $G$ by adding a single free vertex then
$Z_{G_1}=fZ_G$. 
A vertex $i\in V_0$ without any bonds can occupy only one state and doesn't
contribute to the energy at all. So it can safely be omitted.

Now fix a boundary bond $b\in B_0$ (connecting site $i_b$
to the wall) and let $G_D$ be the graph obtained
from $G$ by deleting the bond and denote by $G_C$ the graph 
where the bond is contracted such that its site is located on 
the boundary. In $G_C$  $i_b$  hence becomes a vertex in $V_0$.

\begin{equation} Z_G=CZ_{G_D}+DZ_{G_C} \label{g0}
\end{equation}

We shall give the short calculation establishing this relation.
The proof of (\ref{g1}) is essentially the same.  
We write the set of states of $G$ as a disjoint union
${\cal S}={\cal S}_b\cup{\cal S}_{\overline{b}}$ with
${\cal S}_b:=\{S\in{\cal S}\mid S_{i_b}=0\}$.
\begin{eqnarray*}
Z_G&=&\sum_{S\in{\cal S}}\prod_{(i,j)\in B_1}(1+\delta(S_i,S_j)B)
  \prod_{i\in B_0\backslash \{i_b\}}(C+\delta(0,S_i)D)(C+D\delta(0,S_{i_b}))\\
&=&\sum_{S\in{\cal S}_b}\prod_{(i,j)\in B_1}(1+\delta(S_i,S_j)B)
  \prod_{i\in B_0\backslash \{i_b\}}(C+\delta(0,S_i)D)(C+D)\\
&&+\sum_{S\in{\cal S}_{\overline{b}}}\prod_{(i,j)\in B_1}(1+\delta(S_i,S_j)B)
  \prod_{i\in B_0\backslash \{i_b\}}(C+\delta(0,S_i)D)C
\end{eqnarray*}
By definition $\sum_{S\in{\cal S}_b}\prod_{(i,j)\in B_1}(1+\delta(S_i,S_j)B)
  \prod_{i\in B_0\backslash\{i_b\}}(C+\delta(0,S_i)D)$
is nothing but $Z_{G_C}$.  
The second summand in the calculation of $Z_G$ is almost $CZ_{G_D}$. 
The difference is due to states that assign $0$ to $S_{i_b}$. However 
the contribution of this states to $Z_{G_D}$ is exactly the same 
as $Z_{G_C}$ so that we can complete our calculation
$Z_G=(C+D)Z_{G_C}+C(Z_{G_D}-Z_{G_C})=DZ_{G_C}+CZ_{G_D}$.
 
The relations shown so far suffice to calculate the partition function for
any boundary graph because we can break any bond and thus
work towards trivial
components. 

Figure \ref{latticeskein} displays these relations graphically. 

\setlength{\unitlength}{0.1mm}

\begin{figure}[ht]
\begin{picture}(1500,600)
\linethickness{0.5mm}

\put(0,500){\mbox{$Z${\huge ( }}}
\put(98,520){\circle*{20}} \put(140,520){\circle*{20}}
\put(68,520){\line(1,0){110}}
\put(140,520){\line(1,1){38}} 
\put(140,520){\line(1,-1){38}}
\put(98,520){\line(-1,1){38}} 
\put(98,520){\line(-1,-1){38}}
\put(410,520){\circle*{20}} \put(450,520){\circle*{20}}
\put(180,500){\mbox{{\huge ) }$= AZ${\huge (} }}
\put(450,520){\line(1,1){38}} 
\put(450,520){\line(1,-1){38}}
\put(450,520){\line(1,0){38}}
\put(410,520){\line(-1,1){38}} 
\put(410,520){\line(-1,-1){38}}
\put(410,520){\line(-1,0){38}}
\put(500,500){\mbox{{\huge )}$+BZ${\huge (}}}
\put(720,520){\circle*{20}}
\put(720,520){\line(1,1){38}} 
\put(720,520){\line(1,-1){38}}
\put(690,520){\line(1,0){74}}
\put(720,520){\line(-1,1){38}} 
\put(720,520){\line(-1,-1){38}}
\put(720,520){\line(-1,0){38}}
\put(770,500){\mbox{{\huge ) }}}

\put(900,500){\mbox{Z{\huge (}}}
\put(980,520){\circle*{20}} 
\put(1020,500){\mbox{G {\huge )}$=fZ(G)$}}


\put(0,300){\mbox{$Z${\huge ( }}} 

\linethickness{1mm}\put(100,280){\line(0,1){100}}\linethickness{0.5mm} 
\multiput(65,272)(0,18){5}{\line(1,1){38}}

\put(150,320){\circle*{20}}
\multiput(150,320)(-10,0){5}{\circle*{4}}
\put(150,320){\line(1,0){38}}
\put(145,320){\line(1,1){38}} 
\put(145,320){\line(1,-1){38}}
\put(450,320){\circle*{20}}
\put(180,300){\mbox{{\huge ) }$= CZ${\huge (} }}
\linethickness{1mm}\put(410,280){\line(0,1){100}}\linethickness{0.5mm} 
\multiput(375,272)(0,18){5}{\line(1,1){38}}
\put(450,320){\line(1,1){38}} 
\put(450,320){\line(1,-1){38}}
\put(450,320){\line(1,0){35}}
\put(500,300){\mbox{{\huge )}+$DZ${\huge (}}}
\put(720,320){\circle*{20}}
\put(720,320){\line(1,1){38}} 
\put(720,320){\line(1,-1){38}}
\put(720,320){\line(1,0){38}}
\linethickness{1mm}\put(720,280){\line(0,1){100}}\linethickness{0.5mm} 
\multiput(685,272)(0,18){5}{\line(1,1){38}}
\put(770,300){\mbox{{\huge ) }}}

\put(900,300){\mbox{$Z${\huge ( }}}
\put(1000,320){\circle*{20}}
\linethickness{1mm}
\put(1000,280){\line(0,1){100}}\linethickness{0.5mm} 
\multiput(965,272)(0,18){5}{\line(1,1){38}}
\put(1020,300){\mbox{$G$ {\huge )}$=Z(G)$}}

\put(0,100){\mbox{$Z${\huge ( }}} 
\linethickness{1mm}\put(100,80){\line(0,1){100}}\linethickness{0.5mm} 
\multiput(65,72)(0,18){5}{\line(1,1){38}}
\put(100,100){\circle*{20}}
\put(100,130){\oval(60,60)[r]}
\put(140,100){\mbox{{\huge ) }{\small $= C+D$}}}

\put(420,100){\mbox{$Z${\huge ( }}} 
\linethickness{1mm}\put(520,80){\line(0,1){100}}\linethickness{0.5mm} 
\multiput(485,72)(0,18){5}{\line(1,1){38}}
\put(520,100){\circle*{20}}
\put(520,160){\circle*{20}}
\put(520,130){\oval(60,60)[r]}
\put(550,100){\mbox{{\huge ) }{\small $= A+B$}}}

\put(800,100){\mbox{$Z${\huge ( }}} 
\linethickness{1mm}\put(900,80){\line(0,1){110}}\linethickness{0.5mm} 
\multiput(865,72)(0,18){5}{\line(1,1){38}}
\put(900,100){\circle*{20}}
\put(940,140){\circle*{20}}
\put(900,100){\line(1,1){40}}
\multiput(900,180)(10,-10){4}{\circle*{5}}
\put(980,100){\mbox{{\huge ) }{\small $=A(D+fC)+$}}}
\put(1090,50){\mbox{$B(C+D)$}}
\end{picture}
\caption{\label{latticeskein} Relations between partition functions of models 
on locally changed graphs}
\end{figure}

\setlength{\unitlength}{1mm}

\section{Links and Temperley-Lieb Algebras of Coxeter Type B}

To proceed we need some facts about type B knot theory.
To every root system (of a simple Lie algebra)
there exists the associated Weyl group (Coxeter group).
For root systems of type $A_n$ it is the permutation group. For type $B_n$
it is a semidirect product of a permutation group with $\ZZ_2^n$.
It has generators $X_0,X_1,\ldots,X_{n-1}$ and relations
$X_i^2=1,|i-j|>1\Rightarrow X_iX_j=X_jX_i,
i+1=j>0\Rightarrow X_jX_iX_j=X_iX_jX_i$
and
$X_0X_1X_0X_1=X_1X_0X_1X_0$.
Omitting the quadratic relations from the Coxeter
presentations of these groups
one obtains the braid group of the root system. T. tom Dieck
initiated in \cite{tD1} the systematic study of
 quotients of the group algebras of these braid groups.
\begin{de} The braid group  ${\rm ZB}_n$ of Coxeter type $B$ 
is generated by $X_0,X_1,\ldots,X_{n-1}$ and relations
\begin{eqnarray}
X_0X_1X_0X_1&=&X_1X_0X_1X_0\\
X_iX_jX_i&=&X_jX_iX_j\qquad |i-j|=1,i,j\geq 1\\
X_iX_j&=&X_jX_i\qquad |i-j|>1 \label{taukomm}
\end{eqnarray}
\end{de}

We further need the group $\widetilde{\ZB_n}$ which is 
the free group on generators 
$\sigma_0,\sigma_1,\ldots,\sigma_{n-1},\sigma'_0,
\sigma'_1,\ldots,\sigma'_{n-1}$
with 
(\ref{taukomm}) 
valid 
with 
$X_i$ 
replaced 
by $\sigma_i$ or
$\sigma'_i$. 
Obviously $\sigma_i\mapsto X_i,\sigma'_i\mapsto X_i^{-1}$ is
a surjection.  

Among the finite dimensional quotients of the group algebra of 
$\ZB_n$ there is the
Temperley Lieb algebra of type $B$ that was studied in \cite{tD1}.
Here ee use  a slight generalisation that depends on more
parameters. 
\begin{de} The Temperley-Lieb Algebra ${\rm TB}_n$ of Coxeter type $B$ 
over a ring with
parameters $c,c',d$ is generated by $e_0,e_1,\ldots,e_{n-1}$ and relations
\begin{eqnarray}
e_1e_0e_1&=&c'e_1\\
e_ie_je_i&=&e_i\qquad |i-j|=1,i,j\geq 1\\
e_ie_j&=&e_je_i\qquad |i-j|>1\\
e_0^2&=&ce_0\\
e_i^2&=&de_i
\end{eqnarray}
\end{de}

${\rm TB}_n$ has a Markov trace recursively defined by:
\begin{equation}
{\rm tr}(1)=1\qquad{\rm tr}(e_0)=cd^{-1}\qquad
{\rm tr}(ae_nb):=d^{-1}{\rm tr}(ab)\quad\forall a,b\in{\rm TB}_n
\end{equation}

There are two  topological  interpretations of these algebras.
The first approach is to interpret the words in this algebras
as symmetric tangles. A braid (or tangle) of type $B$  in ${\rm ZB}_n$ has 
$2n$ strands starting and ending in two sets of points (each
set  numbered $\{-n,\ldots,-1,1,\ldots,n\}$) and is symmetric under 
reflection at the middle plane.
  
\begin{figure}[ht]
\begin{picture}(150,60)
\put(72,50){\mbox{$X_0$}}

\linethickness{0.2mm}
\put(2,50){\mbox{$\cdots$}}
\put(9,56){\mbox{{\small -3}}}
\put(14,56){\mbox{{\small -2}}}
\put(19,56){\mbox{{\small -1}}}
\put(40,56){\mbox{{\small 3}}}
\put(35,56){\mbox{{\small 2}}}
\put(30,56){\mbox{{\small 1}}}

\put(10,45){\line(0,1){10}}
\put(15,45){\line(0,1){10}}
\put(20,55){\line(1,-1){10}}
\put(20,45){\line(1,1){4}}
\put(26,51){\line(1,1){4}}
\put(35,45){\line(0,1){10}}
\put(40,45){\line(0,1){10}}
\put(42,50){\mbox{$\cdots$}}

\linethickness{0.4mm}
\put(90,45){\line(0,1){3}}
\put(90,50){\line(0,1){5}}
\linethickness{0.2mm}
\put(88,51){\oval(4,4)[l]}
\put(88,49){\line(1,0){5}}
\put(93,47){\oval(4,4)[tr]}
\put(93,55){\oval(4,4)[br]}
\put(95,56){\mbox{{\small 1}}}
\put(90,56){\mbox{{\small 0}}}

\put(100,45){\line(0,1){10}}

\put(107,50){\mbox{$\cdots$}}

\put(72,35){\mbox{$X_i$}}

\put(0,40){\line(1,-1){10}}
\put(0,30){\line(1,1){4}}
\put(6,36){\line(1,1){4}}
\put(45,40){\line(1,-1){10}}
\put(45,30){\line(1,1){4}}
\put(51,36){\line(1,1){4}}
\put(20,30){\line(0,1){10}}
\put(30,30){\line(0,1){10}}
\put(13,35){\mbox{$\cdots$}}
\put(35,35){\mbox{$\cdots$}}
\linethickness{0.4mm}
\put(90,30){\line(0,1){10}}
\linethickness{0.2mm}
\put(95,30){\line(0,1){10}}
\put(98,35){\mbox{$\cdots$}}
\put(102,40){\line(1,-1){10}}
\put(102,30){\line(1,1){4}}
\put(108,36){\line(1,1){4}}

\put(72,20){\mbox{$e_0$}}

\linethickness{0.2mm}
\put(2,20){\mbox{$\cdots$}}
\put(9,26){\mbox{{\small -3}}}
\put(14,26){\mbox{{\small -2}}}
\put(19,26){\mbox{{\small -1}}}
\put(40,26){\mbox{{\small 3}}}
\put(35,26){\mbox{{\small 2}}}
\put(30,26){\mbox{{\small 1}}}

\put(10,15){\line(0,1){10}}
\put(15,15){\line(0,1){10}}
\put(25,15){\oval(8,8)[t]}
\put(25,25){\oval(8,8)[b]}
\put(35,15){\line(0,1){10}}
\put(40,15){\line(0,1){10}}
\put(42,20){\mbox{$\cdots$}}

\linethickness{0.4mm}
\put(90,15){\line(0,1){10}}
\linethickness{0.2mm}
\put(90,17){\oval(4,4)[tr]}
\put(90,25){\oval(4,4)[br]}
\put(93,26){\mbox{{\small 1}}}
\put(89,26){\mbox{{\small 0}}}

\put(100,15){\line(0,1){10}}

\put(107,20){\mbox{$\cdots$}}

\put(72,5){\mbox{$e_i$}}
\linethickness{0.4mm}
\put(90,0){\line(0,1){10}}
\linethickness{0.2mm}
\put(95,0){\line(0,1){10}}
\put(98,5){\mbox{$\cdots$}}
\put(106,0){\oval(8,8)[t]}
\put(106,10){\oval(8,8)[b]}
\put(20,0){\line(0,1){10}}
\put(30,0){\line(0,1){10}}
\put(12,5){\mbox{$\cdots$}}
\put(32,5){\mbox{$\cdots$}}
\put(5,0){\oval(8,8)[t]}
\put(5,10){\oval(8,8)[b]}
\put(45,0){\oval(8,8)[t]}
\put(45,10){\oval(8,8)[b]}

\end{picture}
\caption{\label{generat} The algebra generators and their interpretation
as symmetric tangles (lefthand side) and as cylinder tangles
 (righthand side) }
\end{figure}
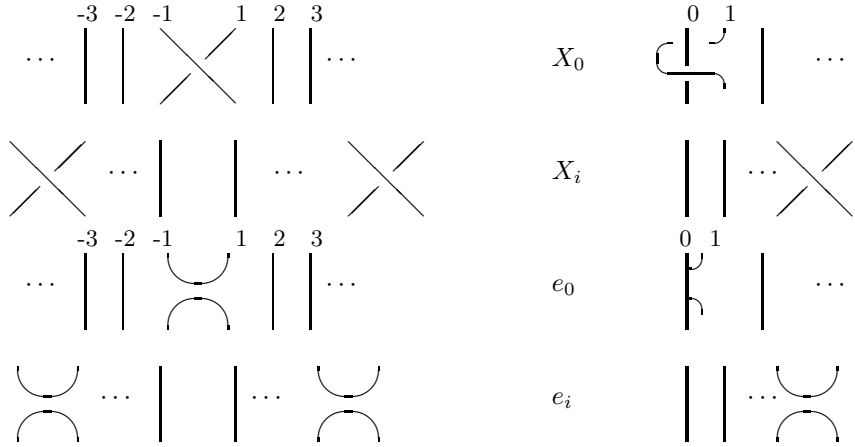

  Alternatively one may think of $B$ type tangles  as
tangles which live in $\RR^3$ where a fixed line (we call it
the 0-strand) has been
removed. We call this second picture the cylinder interpretation.
Both interpretations are useful in our context and
we illustrate them in  Figure \ref{generat}.

From such a purely topological point of view it is natural to require
$c=c'$
in the definition of $\TB_n$. Keeping $c$ and $c'$ 
independent means that we keep information
about the angle between strands $0$ and $1$. 

The map $\phi:\widetilde{{\rm ZB}_n}\rightarrow{\rm TB}_n$ defined by
\begin{eqnarray} \phi(\sigma_i)&:=&\beta e_i+\alpha\qquad
                 \phi(\sigma'_i):=\alpha e_i+\beta\qquad i\geq1\\
\phi(\sigma_0)&:=&\beta_0 e_0+\alpha_0\qquad
                 \phi(\sigma'_0):=\alpha_0 e_0+\beta_0
\end{eqnarray}
 can be graphically
interpreted as implementation of the  skein relations
shown in figure \ref{skein} (using the cylinder picture).

\begin{figure}[ht]
\begin{picture}(150,45)
\linethickness{0.2mm}
\put(15,40){\line(1,-1){10}}
\put(15,30){\line(1,1){4}}
\put(21,36){\line(1,1){4}}
\put(30,35){\mbox{$=\qquad\beta$}}
\put(55,30){\oval(8,8)[t]}
\put(55,40){\oval(8,8)[b]}
\put(70,35){\mbox{$+\alpha$}}
\put(77,30){\line(0,1){10}}
\put(82,30){\line(0,1){10}}

\linethickness{0.4mm}
\put(20,15){\line(0,1){3}}
\put(20,20){\line(0,1){5}}
\linethickness{0.2mm}
\put(18,21){\oval(4,4)[l]}
\put(18,19){\line(1,0){5}}
\put(18,23){\line(1,0){5}}
\put(23,17){\oval(4,4)[tr]}
\put(23,25){\oval(4,4)[br]}

\put(30,20){\mbox{$=\qquad \beta_0$}}
\linethickness{0.4mm}\put(57,15){\line(0,1){10}}\linethickness{0.2mm}
\put(57,15){\oval(8,8)[tr]}
\put(57,25){\oval(8,8)[br]}

\put(80,20){\mbox{$+\alpha_0$}}
\linethickness{0.4mm}\put(90,15){\line(0,1){10}}\linethickness{0.2mm}
\put(95,15){\line(0,1){10}}

\linethickness{0.4mm}\put(20,0){\line(0,1){10}}\linethickness{0.2mm}
\qbezier(20,2)(24,5)(20,8)
\put(27,5){\mbox{=c}}
\linethickness{0.4mm}\put(36,0){\line(0,1){10}}\linethickness{0.2mm}

\linethickness{0.4mm}\put(50,0){\line(0,1){10}}\linethickness{0.2mm}
\put(50,2){\oval(2,2)[tr]}\put(51.5,2){\oval(1,1)[b]}
\put(50,8){\oval(2,2)[br]}\put(51.5,8){\oval(1,1)[t]}
\put(52,2){\line(0,1){6}}
\put(57,5){\mbox{=c'}} 
\linethickness{0.4mm}\put(66,0){\line(0,1){10}}\linethickness{0.2mm}

\put(92,5){\oval(8,8)}\put(100,5){\mbox{$=d$}}

\end{picture}
\caption{\label{skein} Skein relations}
\end{figure}

\section{The Link of a Graph with Boundary}

In analogy with Kauffman's treatment of the ordinary Potts
 model \cite{kauff}
we associate a link diagram $L(G)$ of type B with the boundary graph $G$. 
To accomplish this we
put a crossing on every bond (internal as well as
boundary bonds) and connect their ends so that
each of the cells we produce
contains either a vertex
or a point of the reflecting wall in which a
boundary bond ends. Now we turn each crossing
on an internal bond to an over crossing or under crossing
according to the following rule: Colour each cell containing a vertex
or a boundary bond ending point in black while the others remain
white. If one traverses the crossing on one of its arcs in such a way
that the black region is first at the righthand side
 and after the crossing on the lefthand side then
 the arc lies on top of the other arc.
 The crossings on boundary bonds are turned into the new
kind of braiding that is represented by the picture for
$X_0$ in the cylinder picture. 
The 0-strand is considered to be the
wall.
If one traverses the crossing on a boundary bond
on one of its arcs in such a way
that the black region is first at the righthand side
 and after the crossing on the lefthand side then
the arc first over crosses the fixed 0-strand and then
under crosses it to return.
Vertices on the wall are represented by introducing 
the picture for the $e_0$ generator. 

Figure \ref{assoc} shows a simple example of this process.

It is important to recall that we have assigned link diagrams 
to boundary graphs. They may be deformed by isotopy, but no
Reidemeister moves are allowed. 
 This is the reason for the introduction of the
group $\widetilde{\ZB_n}$.   

\begin{figure}[ht]
\begin{picture}(150,30)

\linethickness{1mm}
\put(10,0){\line(0,1){16}}
\put(20,5){\circle*{2}}
\put(30,5){\circle*{2}}
\put(30,15){\circle*{2}}
\linethickness{0.5mm}
\put(30,5){\line(0,1){10}}
\put(20,5){\line(1,0){10}}
\put(20,5){\line(1,1){10}}
\dottedline{10}{5}{1}{0}{10}

\linethickness{1mm}
\put(45,0){\line(0,1){16}}
\put(55,5){\circle*{2}}
\put(65,5){\circle*{2}}
\put(65,15){\circle*{2}}
\linethickness{0.5mm}
\put(55,5){\line(1,0){10}}
\put(65,5){\line(0,1){10}}
\put(55,5){\line(1,1){10}}
\dottedline{45}{5}{1}{0}{10}
\linethickness{0.2mm}
\put(65,5){\oval(10,5)[br]}
\put(65,7.5){\oval(10,5)[tr]}
\put(47.5,5){\oval(8,5)[l]}
\put(47.5,2.5){\line(1,1){7.5}}
\put(47.5,7.5){\line(1,-1){7.5}}
\put(65,15){\oval(10,5)[t]}
\put(55,10){\line(1,0){12.5}}
\put(55,0){\line(1,1){15}}
\put(60,15){\line(0,-1){10}}\put(65,5){\oval(10,5)[bl]}
\put(70,5){\line(0,1){2.5}}

\linethickness{1mm}
\put(85,0){\line(0,1){2}}
\put(85,4){\line(0,1){12}}
\linethickness{0.2mm}
\put(105,5){\oval(10,5)[br]}
\put(105,7.5){\oval(10,5)[tr]}
\put(87.5,5){\oval(8,5)[l]}
\put(105,15){\oval(10,5)[t]}
\put(95,10){\line(1,0){12.5}}
\put(87.5,2.5){\line(1,0){12}}\put(97.5,2.5){\line(1,1){2}}
\put(87.5,7.5){\line(2,1){5}}\put(92.5,10){\line(1,0){2.5}}
\put(100,15){\line(0,-1){4}}\put(100,9){\line(0,-1){4}}
\put(105,5){\oval(10,5)[bl]}
\put(110,5){\line(0,1){2.5}}
\put(110,15){\line(0,-1){2.5}}
\put(101,3.5){\line(1,1){6}}
\put(108,10.5){\line(1,1){3}}

\end{picture}
\caption{\label{assoc} Associating a B link diagram to a boundary graph}
\end{figure}
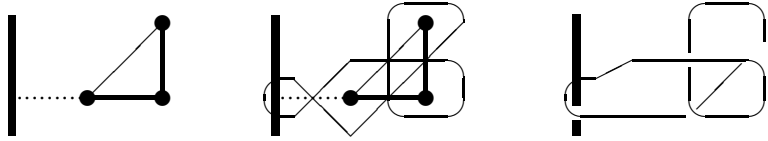

\begin{figure}[ht]
\begin{picture}(150,40)
\linethickness{0.5mm}

\put(20,10){\circle*{2}}
\put(20,20){\circle*{2}}
\put(20,30){\circle*{2}}
\put(30,10){\circle*{2}}
\put(30,20){\circle*{2}}
\put(30,30){\circle*{2}}
\put(40,10){\circle*{2}}
\put(40,20){\circle*{2}}
\put(40,30){\circle*{2}}

\put(20,10){\line(0,1){20}}
\put(30,10){\line(0,1){20}}
\put(40,10){\line(0,1){20}}

\put(20,10){\line(1,0){20}}
\dottedline{10}{10}{1}{0}{10}
\put(20,20){\line(1,0){20}}
\dottedline{10}{20}{1}{0}{10}
\put(20,30){\line(1,0){20}}
\dottedline{10}{30}{1}{0}{10}

\multiput(5,4)(0,3){10}{\line(1,1){5}}

\linethickness{1mm}
\put(10,5){\line(0,1){30}}

\linethickness{0.2mm}

\put(20,35){\line(1,-1){8}}
\put(21,24){\line(1,-1){8}}
\put(21,14){\line(1,-1){9}}
\put(30,35){\line(1,-1){8}}
\put(31,24){\line(1,-1){8}}
\put(31,14){\line(1,-1){9}}

\put(16,11){\line(1,1){8}}
\put(16,21){\line(1,1){8}}
\put(26,11){\line(1,1){8}}
\put(26,21){\line(1,1){8}}
\put(36,11){\line(1,1){9}}
\put(36,21){\line(1,1){9}}

\put(20,5){\line(1,1){4}}
\put(30,5){\line(1,1){4}}

\put(26,31){\line(1,1){4}}
\put(36,31){\line(1,1){4}}

\put(41,14){\line(1,-1){4}}
\put(41,24){\line(1,-1){4}}

\put(40,35){\line(1,-1){5}}

\put(40,5){\line(1,1){5}}

\linethickness{0.4mm}
\put(15,0){\line(0,1){7}}
\put(15,9){\line(0,1){8}}
\put(15,19){\line(0,1){8}}
\put(15,29){\line(0,1){8}}
\linethickness{0.2mm}
\put(14,9.5){\oval(3,3)[l]}
\put(16,8){\line(1,-1){3.5}}
\put(14,19.5){\oval(3,3)[l]}
\put(16,18){\line(1,-1){3.5}}
\put(14,29.5){\oval(3,3)[l]}
\put(16,28){\line(1,-1){3.5}}
\put(20,35){\line(-1,-1){4}}

\end{picture}
\caption{\label{latticelink} A lattice with boundary 
together with its graph}
\end{figure}
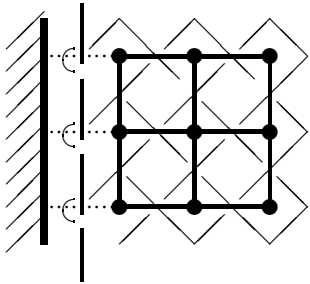

\section{The $B$ Potts Polynomial}

To each link diagram $L$ of type $B$ we associate a 
Polynomial $W(L)$ which is an invariant of isotopy 
(without any Reidemeister moves).
The partition function of a boundary Potts model on
a boundary graph $G$ is then given by
\begin{equation}\label{wpart}
Z_G=C^{\# B_0}d^{\#V}c^{-\#V_0}W(L(G))
\end{equation} 

The $B$ Potts bracket $W(L)$ of a link diagram $L$
is defined by the skein relations of figure \ref{skein}
with the particular values 
\begin{eqnarray}
\alpha&:=&1\qquad\quad\alpha_0:=1\\
\beta&:=&Bd^{-1}\qquad \beta_0:=DC^{-1}c^{-1}\\
d&:=&f^{1/2}\qquad\quad c:=c'd
\end{eqnarray}

The proof of (\ref{wpart}) is done simply by checking
that (\ref{wpart}) changes in the right way (i.e. according to the rules
found in the end of section \ref{pottsmod}) under changes
of the underlying graph. 

Adding a single free vertex in $V\backslash V_0$ enriches
the link diagram by an extra unknotted circle. Thus 
$W$ picks up a factor $d$. An additional factor $d$ is
explicitly multiplied  on the right hand side of (\ref{wpart})
so that $Z$ is multiplied by $d^2=f$ as it should.
If we add an additional vertex $i\in V_0$ on the wall 
the link diagram picks up a loop formed from $e_0$ i.e.
a factor $c'$ according to the skein relations. 
But this factor is cancelled by the contribution $dc^{-1}$ 
arising from the counting of all vertices in (\ref{wpart})
so that $Z$ remains unaltered.

Deleting an inner bond is the same as replacing a 
crossing by the trivial braid. Since this doesn't change
the number of vertices no factor is picked up and since
the skein coefficients for this case agree ($A=\alpha=1$)
the first summand in (\ref{g1}) is reproduced correctly.
Contracting a bond in the graph means introduction of 
parallel strands of the kind of the $e_i$ picture. 
Contracting a bond decreases the number of vertices in 
$V\backslash V_0$ by one. This is taken into account by setting
$\beta:=B/d$.

Now delete a boundary bond. In terms of the link diagram this
means replacing a crossing of type $X_0$ by a trivial braid.
The skein relations give for this a coefficient $\alpha_0=1$. Inserting
the bond back increases the number of boundary bonds by one
so that (\ref{wpart}) introduces an additional factor $C$ as in
(\ref{g0}). For the same reason contracting a boundary bond also
yields a factor $C$ but there is no such factor in the
corresponding summand of (\ref{g0}) hence we have to introduce
a factor $C$ into the  denominator of  
$\beta_0$.  A factor  $D$ in the numerator 
is required by  (\ref{g0}). There is however an additional effect:
When deleting  the bond we increase the number of sites on the wall
by one. This means that the partition function for the original graph
picks up a factor $c^{-1}$ which we also have to supply in $\beta_0$.     

Sine these relations suffice to calculate $Z$ we have shown that 
(\ref{wpart}) indeed expresses the partition function 
considered in section \ref{pottsmod}.

These results can now be applied to the case of a rectangular
lattice. From figure \ref{latticelink} it can be easily seen how
the link diagram of such a lattice can be described in terms 
of the $B$ braid generators.  
\begin{eqnarray}
\tau_{n,m}&:=&\tau'_m(\tau''_m\tau'_m)^n\\
\tau'_m&:=&\sigma_0\sigma_2\cdots \sigma_{2m-2}\\
\tau''_m&:=&\sigma'_1\sigma'_3\cdots \sigma'_{2m-1}\\
E_m&:=&e_1e_3\cdots e_{2m-1}
\end{eqnarray}
$E_m\tau_{n,m}E_m$ is obviously the B link diagram  associated to a
Potts model with boundary living on a $m\times n$ lattice
(with $n$ boundary sites). Hence we have expressed its
$B$ Potts polynomial as
\begin{equation} W=d^{m}{\rm tr}(E_m\phi(\tau_{n,m})E_m)
\end{equation}
and hence its partition function as
\begin{equation} \label{finres}
Z=C^nf^{nm/2}d^{m}{\rm tr}(E_m\phi(\tau_{n,m})E_m)
\end{equation}

\section{Comments}

\begin{enumerate}
\item Knot theory of type $B$ can also be used to find solutions 
of the reflection equation augmenting the parameter dependent Yang Baxter
equations. This was done in \cite{rho} using a Birman Wenzl algebra
of type $B$. Together with the present work this supports our belief
that a great part of the braid and knot program carried out for models on
the plane can also be done for models on the half plane.  

\item It is possible to use (\ref{finres}) to calculate the partition 
function for small lattices. This was done with a Mathematica program
and provided a quick correctness check.
The runtime behaviour of this approach in its direct incarnation is however
worse than that of a direct summation according to (\ref{zdef}). 
\end{enumerate}

\twocolumn

\end{document}